\documentclass{article}

\usepackage[a4paper,top=4.0cm,bottom=2.0cm,left=2.0cm,right=2.0cm]{geometry}

\usepackage{amsmath,amsfonts,mathrsfs,amsfonts}
\usepackage{epsfig,graphicx}

\begin{document}

\title{A hybrid ensemble transform particle filter for nonlinear and spatially extended 
dynamical systems}

\author{Nawinda Chustagulprom\thanks{Universit\"at Potsdam, 
Institut f\"ur Mathematik, Karl-Liebknecht-Str. 24/25, D-14476 Potsdam, Germany} 
\and Sebastian Reich\thanks{Universit\"at Potsdam, 
Institut f\"ur Mathematik, Karl-Liebknecht-Str. 24/25, D-14476 Potsdam, Germany ({\tt sreich@math.uni-potsdam.de}) and University of Reading, Department of Mathematics and Statistics, Whiteknights, PO Box 220, Reading RG6 6AX, UK.} \and Maria Reinhardt \thanks{Universit\"at Potsdam, 
Institut f\"ur Mathematik, Karl-Liebknecht-Str. 24/25, D-14476 Potsdam, Germany}
}

\maketitle

\begin{abstract} Data assimilation is the task to combine evolution 
models and observational data in order to produce
reliable predictions. In this paper, we focus on ensemble-based recursive data assimilation 
problems. Our main contribution
is a hybrid filter  that allows one to adaptively bridge between ensemble Kalman and particle filters. 
While ensemble Kalman filters are robust and applicable to strongly nonlinear systems even with small and moderate ensemble sizes, 
particle filters are asymptotically consistent in the large ensemble size limit. We demonstrate numerically that our hybrid approach can 
improve the performance of both Kalman and particle filters at moderate ensemble sizes. We also show how to implement the concept of 
localization into a hybrid filter, which is key to its applicability to spatially extended systems. 
\end{abstract}


\noindent {\bf AMS(MOS) subject classifications.} 65C05, 62M20, 93E11, 62F15, 86A22


\section{Introduction}

This paper is concerned with algorithms for assimilating observational data into computational 
models using ensemble prediction techniques and sequential Bayesian model updates. 
The focus will be on state estimation problems which give rise to the classic filtering
problem of stochastic processes \cite{sr:jazwinski}. 

There are two popular classes of ensemble based filtering techniques; namely particle
filters (PFs) \cite{sr:Doucet} and ensemble Kalman filters (EnKFs) \cite{sr:evensen}.
While PFs lead to consistent approximations of the underlying filtering problem 
\cite{sr:crisan}, they suffer from the curse of dimensionality 
\cite{sr:bengtsson08}. EnKFs, on the other hand, are applicable
to spatially extended systems when combined with inflation and localization \cite{sr:evensen}; 
but they lead to estimates which are inconsistent with the Bayesian approach to filtering unless all involved probability density
functions (PDFs) are Gaussian. There have been several recent attempts to enhance the
performance of PFs by either using better proposal densities (see \cite{sr:leeuwen15} for a review)
or by implementing localisation into PFs \cite{sr:vanhandel15,sr:reich15}. 
At the same time, several efforts
have been made to extend EnKFs beyond their Gaussian prior assumption. Most of these attempts 
are built around Gaussian mixture approximations and, hence, can be seen as efforts 
to bridge EnKFs with PFs (see, for example, \cite{sr:smith07,sr:stordal11,sr:reich11,sr:frei11}). 
Other recent attempts to extend EnKFs beyond Gausianity include 
\cite{sr:anderson10,sr:lei11,sr:snyder13}.

In this paper, we combine the general linear ensemble transform filter (LETF) framework, as
laid out in detail in \cite{sr:reichcotter15}, with the EnKF-PF bridging approach suggested
in \cite{sr:frei11}.  As in \cite{sr:frei11}, our hybrid ensemble filter approach bridges 
EnKFs with PFs through a parameter $\alpha \in [0,1]$, which can be chosen adaptively.  
Key innovations of our approach, as compared to the hybrid method of \cite{sr:frei11}, 
include the use of $R$-localisation \cite{sr:hunt07,sr:reichcotter15}, 
the fact that either a PF or an EnKF can be applied 
first, that any implementation of an EnKF can be used (including \cite{sr:anderson10}), 
and that observation operators need not to be linear. 

The properties of our hybrid ensemble filter will be demonstrated numerically for the Lorenz-63 model
\cite{sr:lorenz63}, the Lorenz-96 model \cite{sr:lorenz96,sr:lorenz98}, and a coupled
Lorenz-96 wave equation model \cite{sr:br10b}. The numerical results demonstrate 
that our hybrid approach leads to reduced root mean square (RMS) 
errors compared to an EnKF even for
ensemble sizes at which a standard PF fails to track the reference solution. 

We now summarize the basic data assimilation problem which will be considered 
throughout this paper. See
\cite{sr:stuart15,sr:reichcotter15} for further background material on filtering and general data
assimilation algorithms. The basic task of sequential filtering (or sequential data assimilation) 
is to estimate a reference trajectory $z_{\rm ref}(t)\in \mathbb{R}^{N_z}$ of an evolution equation, 
such as
\begin{equation} \label{model}
\frac{{\rm d}z}{{\rm d}t} = f(z,t),
\end{equation}
from partial and noisy observations
\begin{equation} \label{obs}
y_{\rm obs}(t_k) = h(z_{\rm ref}(t_k),t) + \xi_k, \qquad k=1,2,\ldots,
\end{equation}
at discrete times $t_k> 0$. Here $h:\mathbb{R}^{N_z} \times \mathbb{R}_{>0} 
\to \mathbb{R}^{N_y}$ is the observation
operator and  $\xi_k \in \mathbb{R}^{N_y}$, $k\ge 1$, 
denote independent realizations of a Gaussian random variable with mean zero
and covariance matrix $R \in \mathbb{R}^{N_y\times N_y}$. We are also given the PDF, $\pi_Z(z,0)$, for 
the initial conditions $z(0)$ of (\ref{model}) at time zero. 
In order to simplify the subsequent discussion, we
will assume that the observation operator $h$ is linear and time independent, i.e., $h(z,t) = Hz$.

EnKFs and PFs both aim at providing approximations to
the conditional distributions $$\pi_Z(z,t|y_{\rm obs}(t_1),\cdots,y_{\rm obs}(t_k))$$ for $t\ge t_k$
in the form of an ensemble $z_i(t)$, $i=1,\ldots,M$. The best estimate for $z_{\rm ref}(t)$ is
provided by the ensemble mean,
\begin{equation}
z_{\rm est}(t) = \frac{1}{M} \sum_{i=1}^M z_i(t).
\end{equation}
While $z_{\rm est}(t)$ from a PF converges to the true expectation value
\begin{equation}
\bar z(t) = \int_{\mathbb{R}^{N_z}} z\,\pi_Z(z,t|y_{\rm obs}(t_1),\cdots,y_{\rm obs}(t_k))\,{\rm d}z
\end{equation}
as $M\to \infty$, the same does not hold true for EnKFs, in general. 
However, it is often found in practice that EnKFs lead to smaller time-averaged RMS errors,
\begin{equation}\label{rmse}
\mbox{RMSE} \,\,= \frac{1}{T} \int_0^T \sqrt{ \frac{1}{N_z} \|z_{\rm est}(t)-z_{\rm ref}(t)\|^2}\, {\rm d}t,
\end{equation}
for small or moderate ensemble sizes, $M$. See \cite{sr:handbook} for more details.

Any ensemble-based filter algorithm alternates between prediction
steps, where an ensemble is propagated under the dynamics (\ref{model}), and analysis steps,
where a forecast ensemble $z_i^{\rm f}=z_i(t_k)$ is transformed into an analysis ensemble 
$z_i^{\rm a}$ based on the available observations at time $t=t_k$. 
The analysis ensemble then provides the
initial conditions for the next prediction step over the time interval $t\in (t_k,t_{k+1}]$. 
It is the implementation of the analysis step that provides the crucial difference 
between the various filter algorithms. Section \ref{sec:LETF} lays out a general framework
for formulating analysis steps and discusses the ensemble square root filter (ESRF) 
\cite{sr:tippett03} and the ensemble transform particle filter (ETPF) 
\cite{sr:reich15} as two specific examples for EnKFs and PFs,
respectively, which will be used throughout this paper. The novel hybrid filter is based 
on a combination of the ESRF and the ETPF and is introduced in detail in Section 
\ref{sec:hybrid}. Its generalization to spatially extended 
systems using localization is formulated in Section \ref{sec:spatial}. The complete hybrid algorithm 
is summarized in Section \ref{sec:alg}, while the numerical results can be found in 
Section \ref{sec:num}.


\section{Summary of the LETF framework} \label{sec:LETF}

Given a forecast ensemble $z_i^{\rm f}$, $i=1,\ldots,M$, at an observation time $t_k$, 
we consider general  linear transformations \cite{sr:reichcotter15}
\begin{equation} \label{LETF}
z_j^{\rm a} = \sum_{i=1}^M z_i^{\rm f} d_{ij}
\end{equation}
into an analysis ensemble $z_j^{\rm a}$, $j=1,\ldots,M$. The coefficients $d_{ij}$ satisfy
\begin{equation}
\sum_{i=1}^M d_{ij} =1
\end{equation}
and depend  on the forecast ensemble and the observed $y_{\rm obs}=y_{\rm obs}(t_k)$. 
We demonstrate in the next two paragraphs how the ESRF and the ETPF fit into the 
LETF framework. See \cite{sr:reichcotter15} for further details and LETF formulations of
other EnKFs in particular.

\subsection{ESRF formulation}

We introduce the forecast ensemble mean
\begin{equation}
\bar z^{\rm f} = \frac{1}{M} \sum_{i=1}^M z_i^{\rm f},
\end{equation}
the $N_z \times M$ matrix of ensemble anomalies
\begin{equation}
A^{\rm f} = \left[(z_1^{\rm f}-\bar z^{\rm f}) \,\, (z_2^{\rm f}-\bar z^{\rm f}) \,\, \cdots \,\,
(z_M^{\rm f}-\bar z^{\rm f})\right],
\end{equation}
the symmetric $M\times M$ matrix
\begin{equation}
S = \left\{ I + \frac{1}{M-1} (HA^{\rm f})^{\rm T} 
R^{-1} H A^{\rm f}\right\}^{-1/2},
\end{equation}
and the weight vector $\hat w \in \mathbb{R}^M$ with entries
\begin{equation}
\hat w_i = \frac{1}{M} - \frac{1}{M-1} e_i^T S^2 (HA^{\rm f})^{\rm T}R^{-1} (
H\bar z^{\rm f}  -y_{\rm obs}).
\end{equation}
Here $e_i$ denotes the $i$th basis vector in $\mathbb{R}^M$. Also note that the entries 
$s_{ij}$ of the matrix $S$ satisfy
\begin{equation}
\sum_{j=1}^M s_{ij} = \sum_{i=1}^M s_{ij} = 1 \quad \mbox{and} \quad \sum_{i=1}^M \hat w_i = 1.
\end{equation}
Following the presentation in \cite{sr:reichcotter15}, 
the analysis ensemble of an ESRF is now given by
\begin{equation}
z_j^{\rm a} = \sum_{i=1}^M z_i^{\rm f} \hat w_i + \sum_{i=1}^M (z_i^{\rm f}-\bar z^{\rm f})
s_{ij} = \sum_{i=1}^M z_i^{\rm f}d_{ij}^{\rm KF}
\end{equation}
with filter coefficients
\begin{equation}
d_{ij}^{\rm KF} = d_{ij}^{\rm KF}(\{z_l^{\rm f}\},y_{\rm obs}) := \hat w_i - \frac{1}{M} + s_{ij},
\end{equation}
which depend on the forecast ensemble and the observation, $y_{\rm obs}$.
We note that the analysis mean satisfies
\begin{equation}
\bar z^{\rm a} = \sum_{i=1}^M z_i^{\rm f} \hat w_i.
\end{equation}
In case of nonlinear observation operators, one would replace $HA^{\rm f}$ by the 
$N_y \times M$ matrix
\begin{equation}
\left[ (h(z_1^{\rm f},t_k)-\bar h^{\rm f})\,\, (h(z_2^{\rm f},t_k) - \bar h^{\rm f}) \,\,\cdots\,\,
(h(z_M^{\rm f},t_k)-\bar h^{\rm f})
\right]
\end{equation}
and $H\bar z^{\rm f}$ by
\begin{equation}
\bar h^{\rm f} = \frac{1}{M}\sum_{i=1}^M h(z_i^{\rm f},t_k)
\end{equation}
in the formulas above. We finally mention that 
the computational complexity of the ESRF scales like ${\cal O}(M^3)$ in the ensemble size $M$. 

\subsection{ETPF formulation} \label{sec:ETPF}

The  ETPF provides an implementation of PFs 
which replaces the resampling step of a sequential importance resampling (SIR) 
PF \cite{sr:Doucet} by a linear transformation
of the form (\ref{LETF}) \cite{sr:reich13,sr:reich15,sr:reichcotter15}. As for a SIR
PF, we define importance weights
\begin{equation} \label{PFweights}
w_i = \frac{\exp \left( -\frac{1}{2}(Hz_i^{\rm f}-y_{\rm obs})^{\rm T} R^{-1} (
Hz_i^{\rm f}-y_{\rm obs}) \right)}{\sum_{j=1}^M \exp \left(
 -\frac{1}{2}(Hz_j^{\rm f}-y_{\rm obs})^{\rm T} R^{-1} (
Hz_j^{\rm f}-y_{\rm obs}) \right)}.
\end{equation}
The required transformation (\ref{LETF}) is then found by
solving the linear transport problem
\begin{equation} \label{LT}
T^\ast = \arg \min_{T \in \Pi} \sum_{i,j=1}^M t_{ij} \|z_i^{\rm f}-z_j^{\rm f}\|^2
\end{equation}
where $\Pi$ denotes the set of $M\times M$ matrices $T$ with non-negative entries $t_{ij}$
subject to 
\begin{equation}
\sum_{i=1}^M t_{ij} = 1, \quad \mbox{and} \quad \sum_{j=1}^M t_{ij} = w_i M.
\end{equation}
Finally, the ETPF defines the analysis ensemble via
\begin{equation}
z_j^{\rm a} = \sum_{i=1}^M z_i^{\rm f} d^{\rm PF}_{ij} 
\end{equation}
with filter coefficients
\begin{equation}
d_{ij}^{\rm PF} = d_{ij}^{\rm PF}(\{z_l^{\rm f}\},y_{\rm obs}) := t_{ij}^\ast 
\end{equation}
and analysis ensemble mean
\begin{equation}
\bar z^{\rm a} = \frac{1}{M} \sum_{j=1}^M z_j^{\rm a} = \frac{1}{M} \sum_{j=1}^M 
\sum_{i=1}^M z_i^{\rm f} t_{ij}^\ast = \sum_{i=1}^M z_i^{\rm f} w_i .
\end{equation}
The consistency of a single ETPF step with Bayes' formula 
has been demonstrated in \cite{sr:reich13}. However, the ETPF can fail to track the reference solution, $z_{\rm ref}(t)$, if
the effective sample size
\begin{equation}
M_{\rm eff} = \frac{1}{\sum_{i=1}^M w_i^2}
\end{equation}
becomes too small. In the extreme case of $M_{\rm eff} = 1$, all analysis ensemble members, $z_j^{\rm a}$, become equal
to $z_{i^\ast}^{\rm f}$, where $z_{i^\ast}^{\rm f}$ 
is the forecast ensemble member with importance weight $w_{i^\ast} = 1$.
This limitation of the ETPF provides the main motivation for our hybrid ensemble transform particle filter.

The linear transport problem (\ref{LT}) can be solved by a linear programming algorithm
\cite{sr:wright99}. In our numerical experiments we used the {\it FastEMD} algorithm of 
\cite{sr:Pele-iccv2009}, which has a computational complexity of ${\cal O}(M^3\,\log M )$. 
We also note that one dimensional transport problems can be 
solved efficiently by sorting the ensemble members \cite{sr:reichcotter15}.

Nonlinear observation operators can be accommodated by the ETPF by simply replacing $Hz_i^{\rm f}$ by $h(z_i^{\rm f},t_k)$ in the computation of the importance weights (\ref{PFweights}).


\subsection{The hybrid ensemble transform particle filter} \label{sec:hybrid}

Following \cite{sr:frei11}, we split the likelihood function into a product of two factors, i.e.,
\begin{equation}
\pi_Y(y_{\rm obs}|z) \propto
\exp \left( -\frac{\alpha}{2}(Hz-y_{\rm obs})^{\rm T} R^{-1} (
Hz -y_{\rm obs}) \right) \times 
\exp \left( -\frac{1-\alpha}{2}(Hz-y_{\rm obs})^{\rm T} R^{-1} (
Hz-y_{\rm obs}) \right)
\end{equation}
for $\alpha \in [0,1]$. The second factor is treated by the ESRF giving rise to
\begin{equation}
S(\alpha) = \left\{ I + \frac{1-\alpha}{M-1} (HA^{\rm f})^{\rm T} 
R^{-1} H A^{\rm f}\right\}^{-1/2},
\end{equation}
weight vector $\hat w(\alpha) \in \mathbb{R}^M$ with entries
\begin{equation}
\hat w_i(\alpha) = \frac{1}{M} - \frac{1-\alpha}{M-1} e_i^T S(\alpha)^2 (HA^{\rm f})^{\rm T}R^{-1} (
H\bar z^{\rm f}  -y_{\rm obs}),
\end{equation}
and filter coefficients 
\begin{equation}
d_{ij}^{\rm KF}(\alpha,\{z_l^{\rm f}\},y_{\rm obs}) 
:= \hat w_i(\alpha) - \frac{1}{M} + s_{ij}(\alpha).
\end{equation}

The first factor is treated by the ETPF with importance weights
\begin{equation} \label{weights}
w_i(\alpha) = \frac{\exp \left( -\frac{\alpha}{2}(Hz_i^{\rm f}-y_{\rm obs})^{\rm T} R^{-1} (
Hz_i^{\rm f}-y_{\rm obs}) \right)}{\sum_{j=1}^M \exp \left(
 -\frac{\alpha}{2}(Hz_j^{\rm f}-y_{\rm obs})^{\rm T} R^{-1} (
Hz_j^{\rm f}-y_{\rm obs}) \right)}
\end{equation}
and the linear transport problem
\begin{equation}
T^\ast(\alpha) = \arg \min_{T \in \Pi} \sum_{i,j=1}^M t_{ij} \|z_i^{\rm f}-z_j^{\rm f}\|^2,
\end{equation}
where $\Pi$ now denotes the set of $M\times M$ matrices $T$ with $t_{ij}\ge 0$
subject to
\begin{equation}
\sum_{i=1}^M t_{ij} = 1, \quad \mbox{and} \quad \sum_{j=1}^M t_{ij} = w_i(\alpha)\,M.
\end{equation}
The coefficients of the associated ETPF are then given by
\begin{equation}
d_{ij}^{\rm PF}(\alpha,\{z_l^{\rm f}\},y_{\rm obs}) := t_{ij}^\ast (\alpha).
\end{equation}

When implementing the hybrid filter, one can decide between two options; either 
one applies the ETPF or the ESRF first. The two resulting implementations can be
summarized as follows:
\begin{itemize}
\item[(A)] ETPF-ESRF
\begin{align} \label{hybrida1}
z^{\rm h}_j &= \sum_{i=1}^M z_i^{\rm f} d_{ij}^{\rm PF}, \qquad d_{ij}^{\rm PF} := d_{ij}^{\rm PF}(\alpha,\{z_l^{\rm f}\},y_{\rm obs}),\\ \label{hybrida2}
z^{\rm a}_j &= \sum_{i=1}^M z_i^{\rm h} d_{ij}^{\rm KF}, \qquad d_{ij}^{\rm KF} := d_{ij}^{\rm KF}(\alpha,\{z_l^{\rm h}\},y_{\rm obs}).
\end{align}
\item[(B)] ESRF-ETPF 
\begin{align} 
z^{\rm h}_j &= \sum_{i=1}^M z_i^{\rm f} d_{ij}^{\rm KF}, \qquad d_{ij}^{\rm KF} := d_{ij}^{\rm KF}(\alpha,\{z_l^{\rm f}\},y_{\rm obs}),\\
z^{\rm a}_j &= \sum_{i=1}^M z_i^{\rm h} d_{ij}^{\rm PF}, \qquad d_{ij}^{\rm PF} := d_{ij}^{\rm PF}(\alpha,\{z_l^{\rm h}\},y_{\rm obs}).
\end{align}
\end{itemize}
Note that (\ref{hybrida2}) uses the updated ensemble from (\ref{hybrida1}) to determined
the coefficients $d_{ij}^{\rm KF}$. The same dependency applies to the ESRF-ETPF.

The bridging parameter $\alpha$ can either be set to a constant value or selected adaptively \cite{sr:frei11} 
such that the ratio  between the effective sample size,
\begin{equation}
M_{\rm eff}(\alpha) := \frac{1}{\sum_{i=1}^M w_i(\alpha)^2}, 
\end{equation}
in the ETPF part of the hybrid filter and the ensemble size, $M$, 
is equal to a fixed reference value $\theta \in [0,1]$, i.e.
\begin{equation}
\frac{M_{\rm ref}(\alpha)}{M} = \theta .
\end{equation}
Note that $\theta = 1$ implies $\alpha = 0$ and $M_{\rm ref}(\alpha)$ is monotonically 
decreasing in $\alpha$. If $M_{\rm ref}(\alpha) > \theta M$ for all $\alpha \in [0,1]$, then
we set $\alpha = 1$.

\subsection{Ensemble inflation and particle rejuvenation}

Filter implementations with small or moderate ensemble sizes often lead to filter divergence, i.e., 
the filter is no longer able to track the reference solution. Filter divergence can often be traced back
to analysis ensembles becoming underdispersive, i.e., the analysis ensemble spread does not properly represent
Bayesian posterior uncertainties. In order to prevent ensembles from becoming underdispersive, one can apply
multiplicative ensemble inflation to the forecast ensemble prior to an analysis step and/or one can 
apply particle rejuvenation 
to the analysis ensemble (see, for example, \cite{sr:evensen,sr:reichcotter15}). 
We used particle rejuvenation in the form of 
\begin{equation} \label{rejuvenation}
z_j^{\rm a} \to z_j^{\rm a} + \sum_{i=1}^M (z_i^{\rm f} - \bar z^{\rm f}) \frac{\beta \xi_{ij}}{\sqrt{M-1}}
\end{equation}
in our numerical experiments. Here $\beta >0$ denotes the rejuvenation parameter and 
the $\xi_{ij}$'s are i.i.d.~Gaussian random variables with mean
zero and variance one.  We also enforce 
\begin{equation}
 \sum_{j=1}^M \xi_{ij} = 0
 \end{equation}
 in order to preserve the ensemble mean, $\bar z^{\rm a}$, under rejuvenation.
 
 Note that the particle rejuvenation step is the 
 only instance in the complete filter algorithm where randomness enters. 
 This randomness 
 provides a safeguard against the creation of identically analysis ensemble members under an 
 ETPF, as discussed towards the end of Section \ref{sec:ETPF}.


\section{Spatially extended systems and localization} \label{sec:spatial} 

At each data assimilation step, we have to now deal with 
ensembles of spatial fields $z(x)\in \mathbb{R}^{N}$ where $x\in \mathbb{R}^d$ denotes 
the spatial variable of dimension $d$. Those fields are provided on a computational grid with grid points 
$x_k\in \mathbb{R}^d$, $k=1,\ldots,K$. 
Spatially discrete fields are formally collected into a state vector 
\begin{equation}
z = (z(x_1)^{\rm T} \,\,z(x_2)^{\rm T}\,\,\cdots \,\, z(x_K)^{\rm T} )^{\rm T}
\end{equation}
of dimension $N_z = N\times K$. A forecast ensemble 
\begin{equation}
z_i^{\rm f} = (z_i^{\rm f}(x_1)^{\rm T} \,\,z_i^{\rm f}(x_2)^{\rm T}\,\,\cdots \,\, 
z_i^{\rm f} (x_K)^{\rm T} )^{\rm T}\in \mathbb{R}^{N_z},
\end{equation} 
$i=1,\ldots,M$, is then transformed at each grid point, $x_k$, according to
\begin{equation} \label{locupdate}
z_j^{\rm a}(x_k) = \sum_{i=1}^M z_i^{\rm f}(x_k) d_{ij}(x_k),  \qquad d_{ij}(x_k) = d_{ij}(x_k,\{z_i^f\},
y_{\rm obs}),
\end{equation}
which gives rise to the analysis ensemble 
\begin{equation}
z_j^{\rm a} = (z_j^{\rm a}(x_1)^{\rm T} \,\,z_j^{\rm a}(x_2)^{\rm T}\,\,\cdots \,\, 
z_j^{\rm a} (x_K)^{\rm T} )^{\rm T}\in \mathbb{R}^{N_z},
\end{equation} 
$j=1,\ldots,M$. It is crucial that
these local updates are implemented in a manner that maintains the spatial regularity of 
the forecast states. This prevents the use of localization 
in a standard SIR PF.

We assume that observations are taken at a discrete subset of grid points $\mathfrak{x}_q$, 
$q=1,\ldots,Q$, and satisfy the observation  model
\[
y_q = {\cal H} z(\mathfrak{x}_q) + \eta_q
\]
with independent measurement errors $\eta_q \sim {\rm N}(0,r)$ and observation operator
${\cal H}:\mathbb{R}^{N} \to \mathbb{R}^{N_y}$. We assume here, for simplicity, $N_y = 1$ and $r>0$ denotes the 
variance of the measurement error. We introduce the more compact notation
\begin{equation}
y_{\rm obs} = H z + \eta
\end{equation}
with $y_{\rm obs} = (y_1,\ldots,y_Q)^{\rm T} 
\in \mathbb{R}^Q$, $\eta = (\eta_1,\ldots,\eta_Q)^{\rm T} \sim {\rm N}(0,R)$, $R = rI$, and 
$H:\mathbb{R}^{N_z} \to \mathbb{R}^{Q}$ appropriately defined.  
If each grid point is observed then we simply have $Q=K$, 
$\mathfrak{x}_q = x_q$ with $q = 1,\ldots, K$. If every second
grid point is observed then $Q = K/2$ and $\mathfrak{x}_q = x_{2q}$ with $q = 1,\ldots,K/2$ etc.

We now extend the concept of $R$-localization \cite{sr:hunt07,sr:reich15,sr:reichcotter15}
to the hybrid filters introduced in Section
\ref{sec:hybrid}. First we define a localized inverse measurement error covariance 
matrix $R(x_k)^{-1} \in
\mathbb{R}^{Q\times Q}$ with
diagonal entries
\begin{equation} \label{locR}
r_{qq}(x_k)^{-1} = \frac{\rho\left(\frac{\|x_k-\mathfrak{x}_q\|}{R_{\rm loc}}\right)}{r},
\end{equation}
$q=1,\ldots,Q$, for given localization radius $R_{\rm loc} > 0$ and compactly supported 
tempering function $\rho(t)$ with $\rho(0)= 1$ and $\rho(t) = 0$ for $t \ge 2$. We
used
\begin{equation}
\rho(t) = \left\{ \begin{array}{ll} 1 - \frac{5}{3} t^2 + \frac{5}{8} t^3 + \frac{1}{2} t^4 - \frac{1}{4}t^5
& \mbox{for}\,\,t\le 1,\\
-\frac{2}{3} t^{-1} + 4 - 5t + \frac{5}{3} t^2 + \frac{5}{8} t^3 - \frac{1}{2} t^4 + \frac{1}{12} t^5 &
\mbox{for}\,\,1\le t \le 2,\\
0 & \mbox{otherwise}\,\,, \end{array} \right.
\end{equation}
in our numerical experiments \cite{sr:gaspari99}. 
The two hybrid algorithms from Section \ref{sec:hybrid} 
are now adjusted as follows.  The ESRF is replaced by the 
LETKF \cite{sr:hunt07}, which leads to
\begin{equation}
S(x_k,\alpha_k) = \left\{ I + \frac{1-\alpha_k}{M-1} (HA^{\rm f})^{\rm T} 
R(x_k)^{-1} H A^{\rm f}\right\}^{-1/2},
\end{equation}
\begin{equation}
\hat w_i(x_k,\alpha_k) = \frac{1}{M} - \frac{1-\alpha_k}{M-1} e_i^T S(x_k,\alpha_k)^2 
(HA^{\rm f})^{\rm T}R(x_k)^{-1} (
H\bar z^{\rm f}  -y_{\rm obs}),
\end{equation}
and 
\begin{equation} \label{letkf}
d_{ij}^{\rm KF}(x_k,\alpha_k,\{z_i^{\rm f}\},y_{\rm obs}) 
:= \hat w_i(x_k,\alpha_k) - \frac{1}{M} + s_{ij}(x_k,\alpha_k).
\end{equation}
The ETPF filter contribution  is modified at each grid point $x_k$ as proposed in
\cite{sr:reich15,sr:reichcotter15}. Specifically, we define localized weights
\begin{equation}
w_i(x_k,\alpha_k) = \frac{\exp \left( -\frac{\alpha_k}{2}
(Hz_i^{\rm f}-y_{\rm obs})^{\rm T} R(x_k)^{-1} (
Hz_i^{\rm f}-y_{\rm obs}) \right)}{\sum_{j=1}^M \exp \left(
 -\frac{\alpha_k}{2}(Hz_j^{\rm f}-y_{\rm obs})^{\rm T} R(x_k)^{-1} (
Hz_j^{\rm f}-y_{\rm obs}) \right)},
\end{equation}
introduce the linear transport problem
\begin{equation} \label{transport1}
T^\ast  = \arg \min_{T \in \Pi} \sum_{i,j=1}^M t_{ij} \|z_i^{\rm f}(x_k)-z_j^{\rm f} (x_k)\|^2
\end{equation}
at each grid point, $x_k$, where $\Pi$ denotes the set of $M\times M$ matrices $T$ with $t_{ij}\ge 0$ subject to
\begin{equation} \label{transport2}
\sum_{i=1}^M t_{ij} = 1, \quad \mbox{and} \quad \sum_{j=1}^M t_{ij} = w_i(x_k,\alpha_k)\,M,
\end{equation}
and define
\begin{equation}
d_{ij}^{\rm PF}(x_k,\alpha_k,\{z_i^{\rm f}\},y_{\rm obs}) = t_{ij}^\ast .
\end{equation}
We note that the distance term in (\ref{transport1}) can be replaced by other distance 
measures involving, for example,
grid point values in the vicinity of $x_k$. See \cite{sr:reich15,sr:reichcotter15} for more details. The cost functional
(\ref{transport1}) has the advantage that the associated transport problem is easy to solve in case $z(x_k)\in \mathbb{R}$, i.e.~$N=1$. 

Finally, the desired localized hybrid filter schemes can be implemented in the following two
variants:
\begin{itemize}
\item[(A)] ETPF-LETKF
\begin{align} \label{hybrid1}
z^{\rm h}_j(x_k) &= \sum_{i=1}^M z_i^{\rm f}(x_k) d_{ij}^{\rm PF}(x_k), \qquad d_{ij}^{\rm PF}(x_k) := 
d_{ij}^{\rm PF}(x_k,\alpha_k,\{z_l^{\rm f}\},y_{\rm obs}), \qquad \mbox{over all grid points,}\\ 
\label{hybrid2}
z^{\rm a}_j(x_k) &= \sum_{i=1}^M z_i^{\rm h}(x_k) d_{ij}^{\rm KF}(x_k), 
\qquad d_{ij}^{\rm KF}(x_k) := d_{ij}^{\rm KF}(x_k,\alpha_k,\{z_l^{\rm h}\},y_{\rm obs}), \qquad 
\mbox{over all
grid points.}
\end{align}
\item[(B)] LETKF-ETPF 
\begin{align} 
z^{\rm h}_j(x_k) &= \sum_{i=1}^M z_i^{\rm f}(x_k) d_{ij}^{\rm KF}(x_k), \qquad d_{ij}^{\rm KF} (x_k)
:= d_{ij}^{\rm KF}(x_k,\alpha_k,\{z_l^{\rm f}\},y_{\rm obs}), \qquad \mbox{over all grid points,}\\
z^{\rm a}_j(x_k) &= \sum_{i=1}^M z_i^{\rm h}(x_k) 
d_{ij}^{\rm PF}(x_k), \qquad d_{ij}^{\rm PF}(x_k) := d_{ij}^{\rm PF}(x_k,\alpha_k,
\{z_l^{\rm h}\},y_{\rm obs}),
\qquad \mbox{over all grid points.}
\end{align}
\end{itemize}
The bridging parameters, $\alpha_k$, can again be either set to a constant value, $\alpha$, for all
grid points or determined adaptively. In the latter case $\alpha_k$ is
determined by
\begin{equation} \label{ratio}
\frac{M_{\rm ref}(x_k\alpha_k)}{M} = \theta, \quad M_{\rm ref}(x_k,\alpha_k) := 
\frac{1}{\sum_{i=1}^M w_i(x_k,\alpha_k)^2},
\end{equation}
for given effective ensemble size ratio, $\theta \in [0,1]$. 
If $M_{\rm ref}(x_k,\alpha_k) > \theta M$ for all $\alpha_k \in [0,1]$, then
we set $\alpha_k = 1$. Finally, particle rejuvenation in the form of 
(\ref{rejuvenation}) is performed at the end of each filter step.


\section{Summary of hybrid filter algorithm} \label{sec:alg}

We summarize the numerical implementation of the proposed hybrid filter  
as it has been used for the numerical experiments below. We explicitly describe the filter implementation
with localization. An implementation without localization follows from straightforward
modifications. The free parameters of the data assimilation steps are the ensemble size, $M>1$, 
the localization radius, $R_{\rm loc}>0$, in (\ref{locR}), the bridging parameters, $\alpha_k \in [0,1]$, 
which are either fixed to a value $\alpha$ or determined adaptively  through the effective sample 
size ratio, 
$\theta \in [0,1]$, in (\ref{ratio}). Any choice of these parameters together with a forecast ensemble
$z_l^{\rm f}$, $l=1,\ldots,M$, and an observation $y_{\rm obs}(t_k)$ 
determine the ETPF-LETKF update step (\ref{hybrid1})-(\ref{hybrid2})
uniquely.  The linear transport problems, defined by (\ref{transport1})-(\ref{transport2})
are solved using the {\sc Matlab} implementation of {\it FastEMD} \cite{sr:Pele-iccv2009}. 

The update step (\ref{hybrid1})-(\ref{hybrid2}) is followed by particle rejuvenation (\ref{rejuvenation}) with
parameter $\beta >0$. Note that rejuvenation is applied to the whole spatially dependent analysis
fields with one and the same set of random numbers $\xi_{ij}$, i.e.
\begin{equation}
z_j^{\rm a}(x_k) \to z_j^{\rm a}(x_k) + \sum_{i=1}^M (z_i^{\rm f}(x_k)-\bar z^{\rm f}(x_k))
\frac{\beta \xi_{ij}}{\sqrt{M-1}}
\end{equation}
for all grid points $x_k$. 

A complete implementation also requires a time-stepping method for the differential equation
(\ref{model}). We used the implicit midpoint rule with step-size denoted by $\Delta t>0$. 
Observations (\ref{obs}) are generated from a reference trajectory, $\{z_{\rm ref}^n\}_{n \ge 0}$, 
which is computed using 
the same time-stepping method in (\ref{model}) and a randomly
selected initial condition $z_{\rm ref}^0 \sim \pi_Z(z,0)$. The initial ensemble, $z_i^0$, 
$i=1,\ldots,M$, is also drawn from the initial PDF, $\pi_Z(z,0)$. 


\section{Numerical experiments} \label{sec:num}

\begin{figure}
\begin{center}
\includegraphics[width=0.45\textwidth]{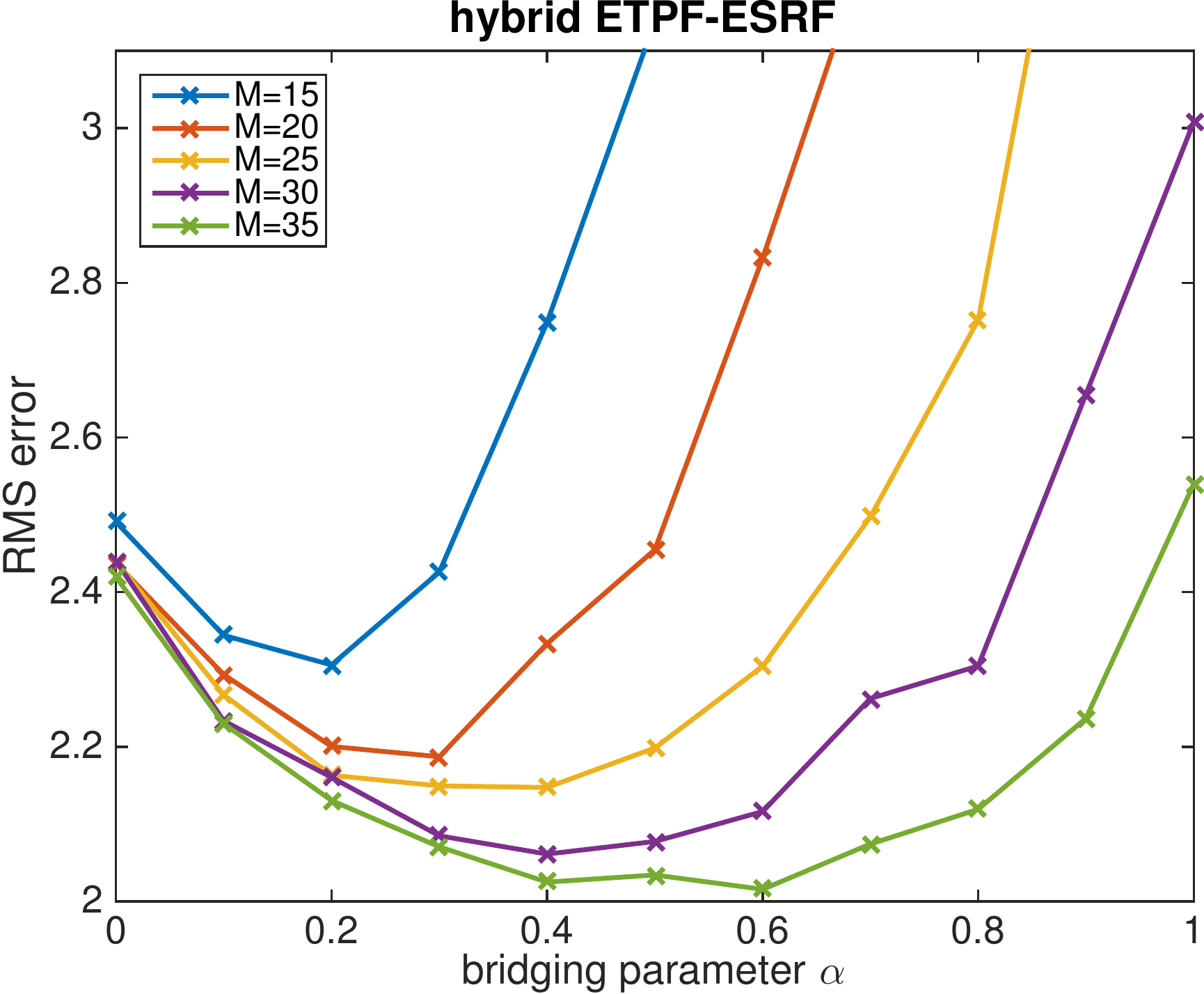} $\quad$
\includegraphics[width=0.45\textwidth]{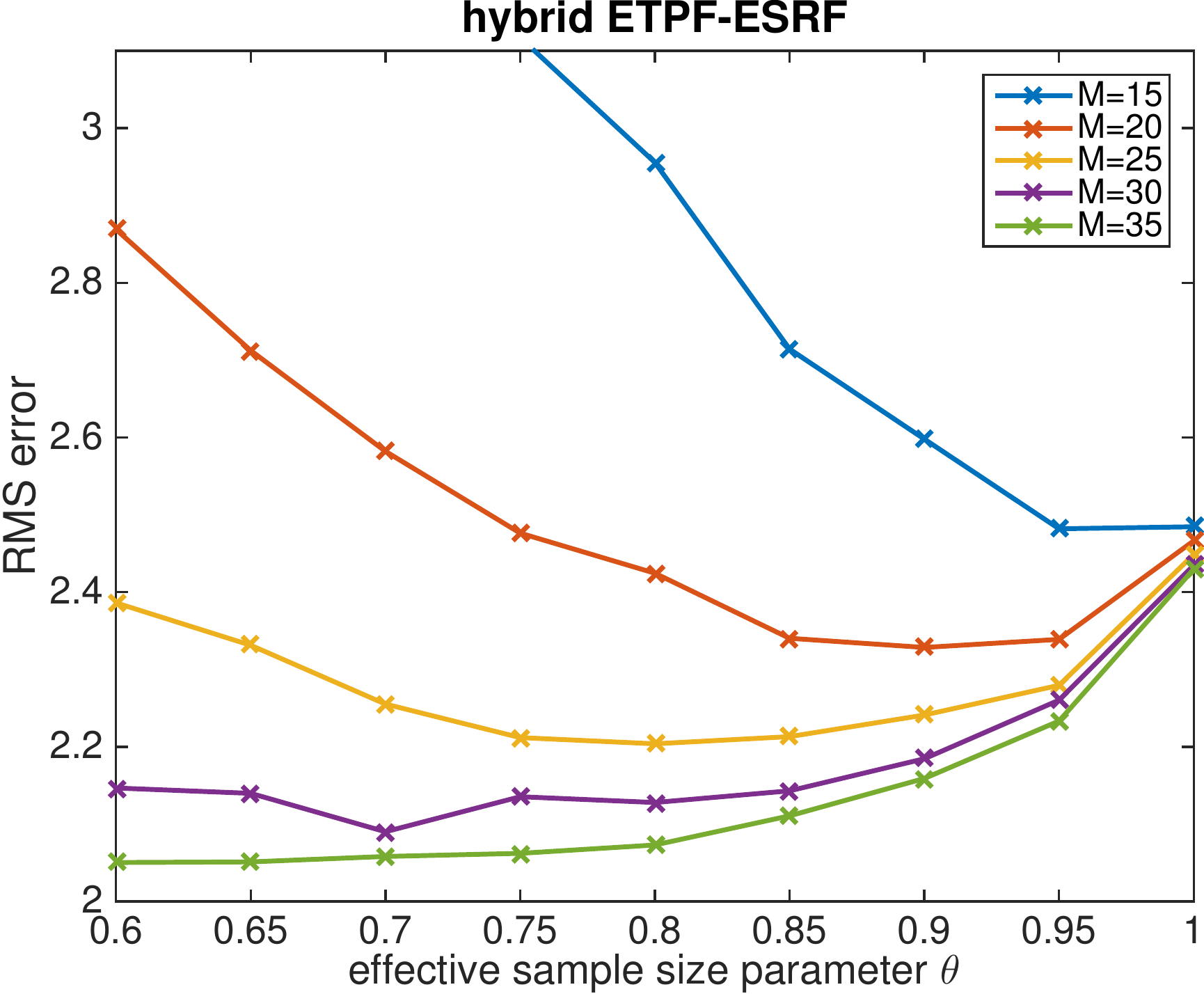}
\end{center}
\caption{Hybrid filter for Lorenz-63 system with ensemble sizes $M=15,\ldots,35$: left panel: RMS
errors for fixed bridging parameter; right panel: RMS errors for adaptively determined 
bridging parameter. Note that
$\alpha = 0$ and $\theta = 1$ correspond to a pure ESRF filter.} \label{figure1}
\end{figure}


We now conduct a series of numerical experiments in order to assess the relative merits of the hybrid 
approach compared to a purely EnKF or PF based data assimilation approach. 
The ensemble sizes, $M$, are chosen
such that standard PF implementations cannot compete with an  EnKF. Such a data assimilation scenario is typical in
complex application areas such as numerical weather prediction \cite{sr:kalnay}.  The quality of a filter is
assessed by the following time-discretised version of (\ref{rmse}):
\begin{equation} \label{RMSE}
\mbox{RMSE} \,\,= \frac{1}{K} \sum_{i=1}^K \sqrt{\frac{1}{N_z} \|\bar z^{\rm a}(t_k)-z_{\rm ref}(t_k)\|^2}.
\end{equation}
Here $K$ denotes the total number of assimilation steps.


\subsection{Lorenz-63}

The first numerical tests of our hybrid filter are performed on the chaotic Lorenz-63 system \cite{sr:lorenz63} with
its standard parameter setting $\sigma = 10$, $\rho = 28$, and $\beta = 8/3$. We only observe the first component
of the three dimensional system in observation intervals of $\Delta t_{\rm obs}=0.12$, i.e., $t_k = 0.12 k$ in 
(\ref{RMSE}), and with observation error variance $R=8$. A total of $K=100,000$ assimilation steps are performed and
$N_z = 3$ in (\ref{RMSE}). The differential equations are solved numerically with a step-size of $\Delta t = 0.01$. 

The ETPF-ESRF filter is implemented for ensemble sizes varying between $M=15$ and $M=35$. Particle rejuvenation
is applied with $\beta = 0.2$. Simulations with $\beta = 0.15$ and $\beta = 0.25$ gave similar results to those reported
here. No localization is necessary for this low dimensional system. Our numerical 
results  for fixed  and adaptively determined bridging parameter, $\alpha$, 
can be found in Figure \ref{figure1}. It can be seen 
that a fixed bridging parameter of $\alpha = 0.2$ leads already to substantial
improvements for the smallest ensemble size of $M=15$. Further improvements are achieved for larger
ensemble sizes at larger bridging parameters. In contrast to these very encouraging results for fixed bridging parameters, the adaptive implementation of the hybrid filter does not work quite as effectively 
for small ensemble sizes of $M=15$ and $M=20$. This finding indicates that one might have
to search for other selection criteria. 

\begin{figure}
\begin{center}
\includegraphics[width=0.45\textwidth]{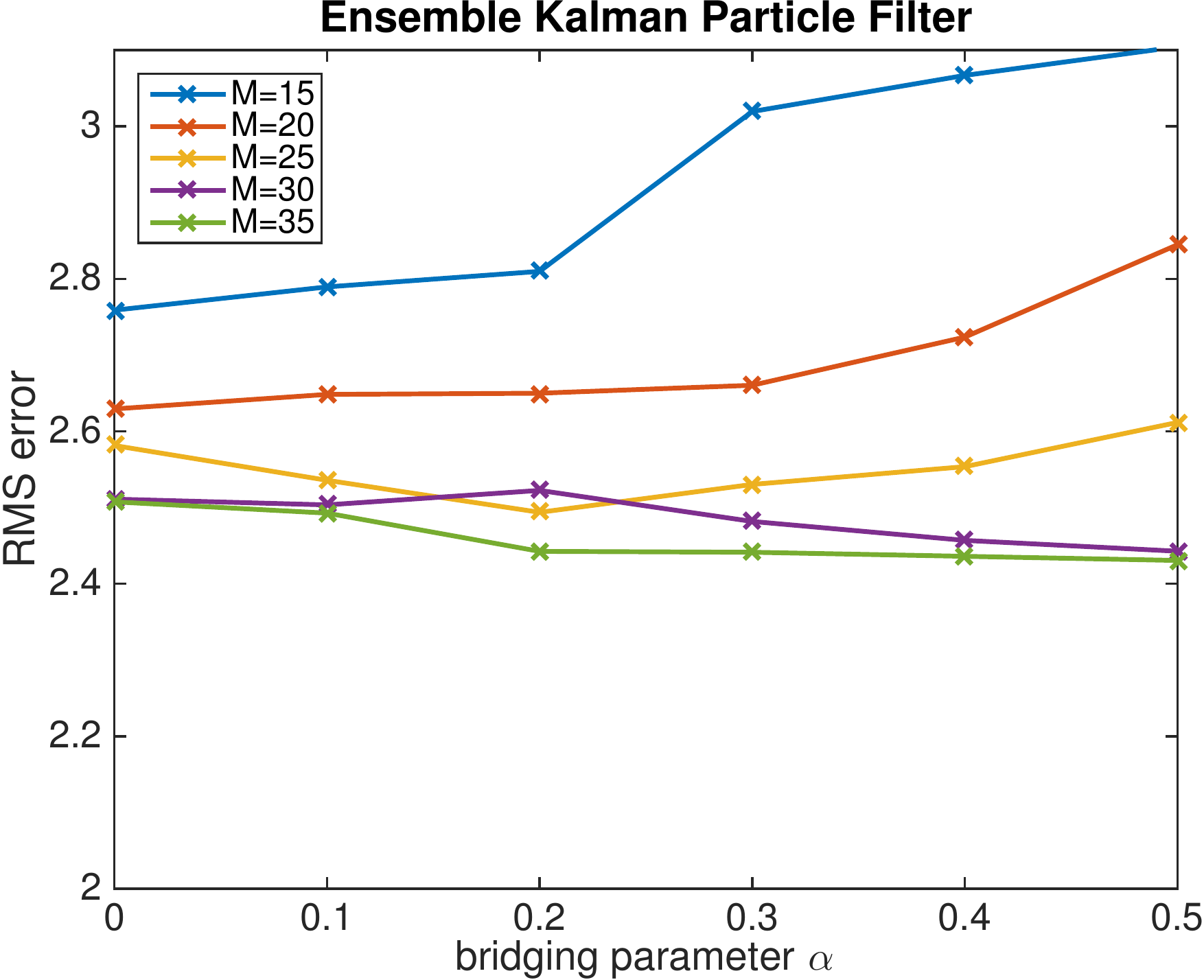} $\quad$
\includegraphics[width=0.45\textwidth]{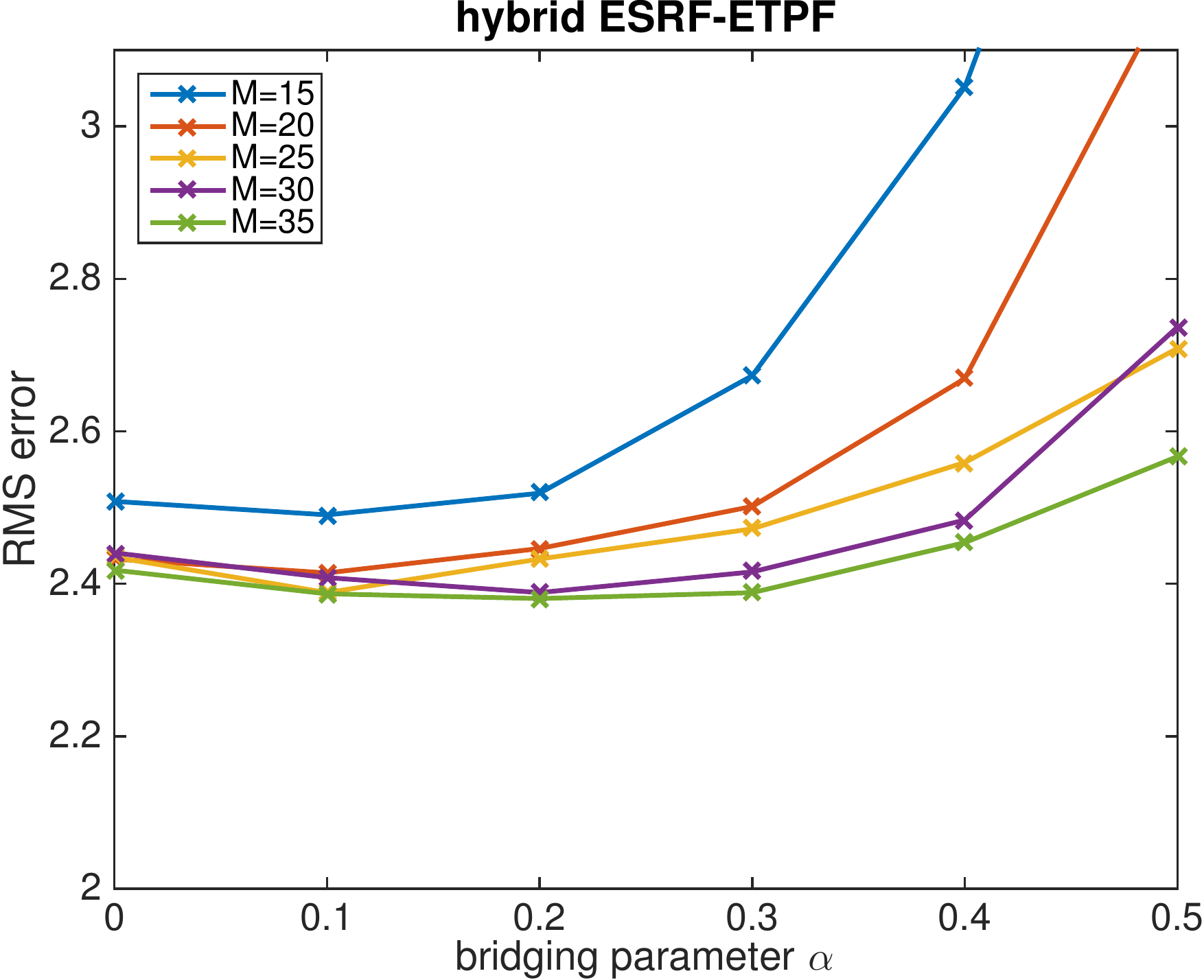}
\end{center}
\caption{Left panel: RMS errors for Ensemble Kalman Particle Filter applied to the 
Lorenz-63 system with ensemble sizes $M=15,\ldots,35$ and for fixed
bridging parameter; right panel: RMS errors for the ESRF-ETPF using the same setting and $\alpha$ in the range
$[0,0.5]$.} \label{figure1b}
\end{figure}

We also implemented the Ensemble Kalman Particle Filter of \cite{sr:frei11} for comparison. The Ensemble Kalman 
Particle Filter bridges the EnKF with perturbed observations \cite{sr:evensen} and a SIR PF. Here we replaced the SIR PF
with the ETPF in order to facilitate a fair comparison. The numerical results are displayed in Figure \ref{figure1b}. 
We  find that the Ensemble Kalman Particle Filter leads to larger RMS errors for all values of the bridging parameter, $\alpha$.
We speculate that this is due to the use of the EnKF with perturbed observations and the fact that the EnKF is applied
first. This hypothesis is supported by the numerical results from the ESRF-ETPF, which also does not lead to significant 
improvements for $\alpha >0$. See Figure \ref{figure1b}. 

Note that we run the hybrid filters at ensemble sizes, $M$, for which 
EnKFs outperform PFs. The more interesting it is that
the ETPF-ESRF leads to substantial improvements for values of the bridging parameter in the range $[0.2,0.4]$ for ensemble sizes
between $M=15$ and $M=35$.  Figure \ref{figure1} shows that the optimal $\alpha$ increases with increasing ensemble sizes, $M$. Asymptotically, one
expects that $\alpha \to 1$ and $\theta \to 0$ as $M\to \infty$. Indeed, we found numerically 
that $\alpha = 0$ is optimal for ensemble sizes smaller than $M=10$ and that $\alpha = 1$ for 
ensemble sizes larger than $M = 100$. See also the convergence study
in Section \ref{sec:convergence}.


\subsection{Lorenz-96}

\begin{figure}
\begin{center}
\includegraphics[width=0.45\textwidth]{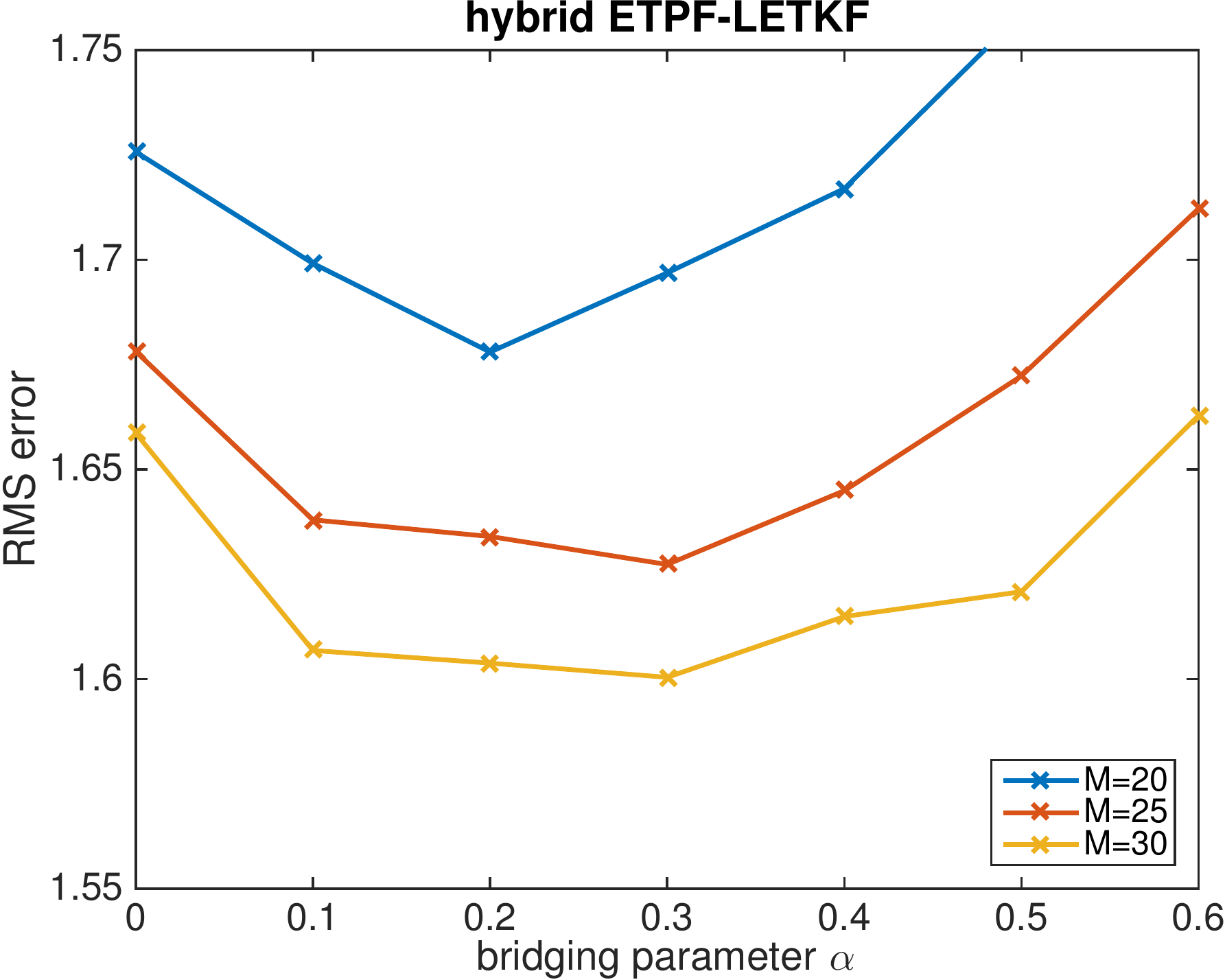} $\quad$
\includegraphics[width=0.45\textwidth]{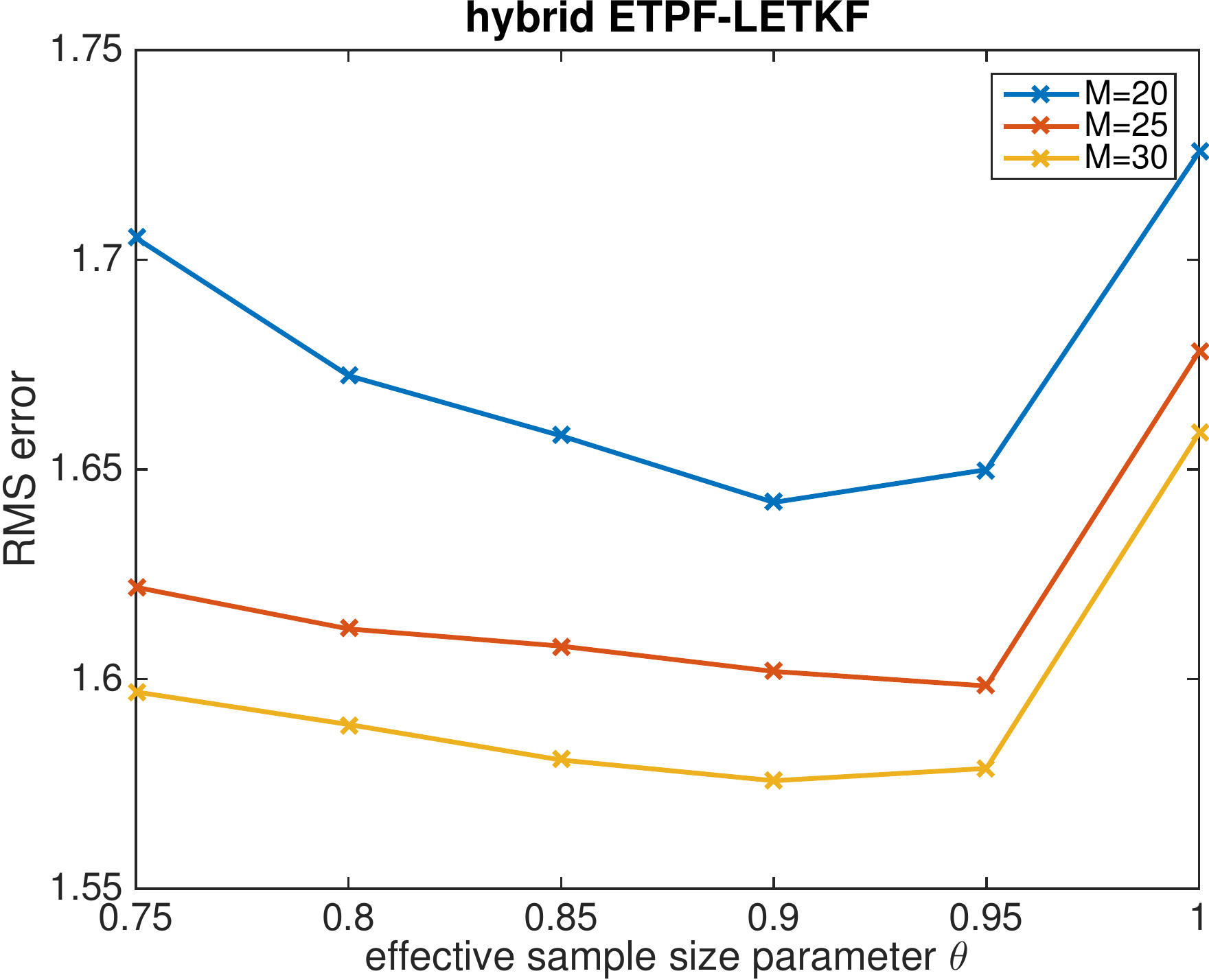}
\end{center}
\caption{Hybrid filter for Lorenz-96 system with ensemble sizes $M=20,25,30$:  left panel: 
RMS errors for fixed bridging parameter; right panel: RMS errors for adaptively determined bridging 
parameters} \label{figure2}
\end{figure}

We next test our hybrid approach on the spatially extended Lorenz-96 model
\begin{align} \label{lorenz96}
\dot{x}_l=(x_{l+1}-x_{l-2})x_{l-1}-x_l+8,\hspace{0.4cm}l=1,...,40
\end{align}
\cite{sr:lorenz96,sr:lorenz98}. 
We observe every other grid point in time intervals of $\Delta t_{\rm obs} = 0.11$ with independent Gaussian 
measurement errors of variance equal to eight. A total of $K=50,000$ assimilation steps is performed. The
hybrid ETPF-LETKF filter is implemented with localization radius $R_{\rm loc} = 4$, particle rejuvenation 
$\beta = 0.2$ and ensemble sizes $M=20,25,30$. The equations of motion are integrated numerically with the
implicit midpoint rule and $\Delta t = \Delta t_{\rm obs}/22$. The resulting time-averaged RMSE can be found in
Figure \ref{figure2}. The overall findings are similar to the Lorenz-63 model except that the adaptive approach
now performs comparably (or even better) to an optimally tuned fixed parameter implementation and that the overall improvements in terms of RMS errors are less significant. The later fact is likely due to the weaker nonlinearity of the Lorenz-96 model. We also implemented the hybrid filter for
localization radii $R_{\rm loc}\in \{3,5,6,7,8\}$ and found that 
$\alpha >0$ always leads to improvements in the RMS error compared to a pure 
LETKF implementation under the simulation conditions stated above.
We finally implemented the Lorenz-96 in the more traditional setting of each grid point being 
observed in intervals of $\Delta t_{\rm obs} = 0.05$ with measurement error variance 
$r = 1$. Our numerical
findings indicate that the achievable improvements of our hybrid filter in comparison to the LETKF
are rather minor and require ensemble sizes $M\ge 30$. We believe that these findings reflect the
fact that frequent observations of the complete state space lead to nearly Gaussian forecast
distributions. 

\subsection{Lorenz-96 coupled to a  wave equation}

Physical processes in the geosciences often act on different time-scales. To simulate that effect we use a model, that consists of the advective Lorenz-96 system (\ref{lorenz96}) and relatively fast moving waves, which are modeled by the discrete wave equation
\begin{align}
\epsilon^2\ddot{h}_l=-h_l+ c^2 (h_{l+1}-2h_l+h_{l-1})-2\epsilon^2\gamma\dot{h}_l,\hspace{0.4cm}l=1,...,40.
\end{align}
The positive parameter $\epsilon \ll 1$ 
characterises the ratio between the advective velocity in $x$ and the wave speed of $h$. The parameter 
$c>0$ determines the wave dispersion over the computational grid and $\gamma>0$ is a damping parameter. The coupled model with coupling-strength $\delta$ looks as follows
\cite{sr:br10b}:
\begin{align}
\dot{x}_l &=(1-\delta)(x_{l+1}-x_{l-2})x_{l-1}+\delta(x_{l-1}h_{l+1}-x_{l-2}h_{l-1})-x_l+8,\\ \label{wave}
\epsilon^2\ddot{h}_l &=-h_l+c^2(h_{l+1}-2h_l+h_{l-1})+x_l-2\epsilon^2\gamma\dot{h}_l,
\end{align}
$l=1,\ldots,40$, and periodic boundary conditions are applied. There exists a balance relation 
\begin{align} \label{balance}
x_l=h_l-c^2(h_{l+1}-2h_l+h_{l-1})
\end{align}
between the $x$ and $h$ variables. If it is fulfilled initially, solutions stay in approximate balance for long time intervals with 
the deviations from exact balance proportional to $\epsilon$. For our experiments, we used this model with a coupling strength of 
$\delta=0.1$, $\epsilon=0.0025$, a damping of $\gamma=0.1$, and $c=1/2$. 
The time-stepping method of \cite{sr:br10b} is used with a step-size of $\Delta t = 0.002$.

The initial ensemble in $\{x_l\}$ is created by sampling from a Gaussian distribution with mean zero and covariance $0.01I$. 
The $h$-values of the initial ensemble are created by calculating them from the $x$-values using the balance relation (\ref{balance}). 
The $\dot{h}$-values are set equal to zero. We observe $x_l$ at every other grid-point at time intervals of 
$\Delta t_{obs} = 0.15$ with measurement-error variance $R = 8 I_{20}$. 
We applied the hybrid ETPF-LETKF filter to our model with a localization radius 
of $R_{\text{loc}}=4$ and a particle-rejuvenation parameter of $\beta = 0.2$.
We let the experiment run for $K=50,000$ assimilation cycles and a spin-up phase of 1000 cycles, in which $\alpha$ is set to zero and the $h$-values are calculated to be in exact balance to the $x$-values after every assimilation step. We conduct experiments 
for ensemble sizes of $M=20,25,30$, fixed bridging parameters $\alpha=0,0.1,...,0.6$ and also for adaptive $\alpha_k$'s using 
$\theta=0.75,0.8,...,1$.

\begin{figure}
\begin{center}
\includegraphics[width=0.45\textwidth]{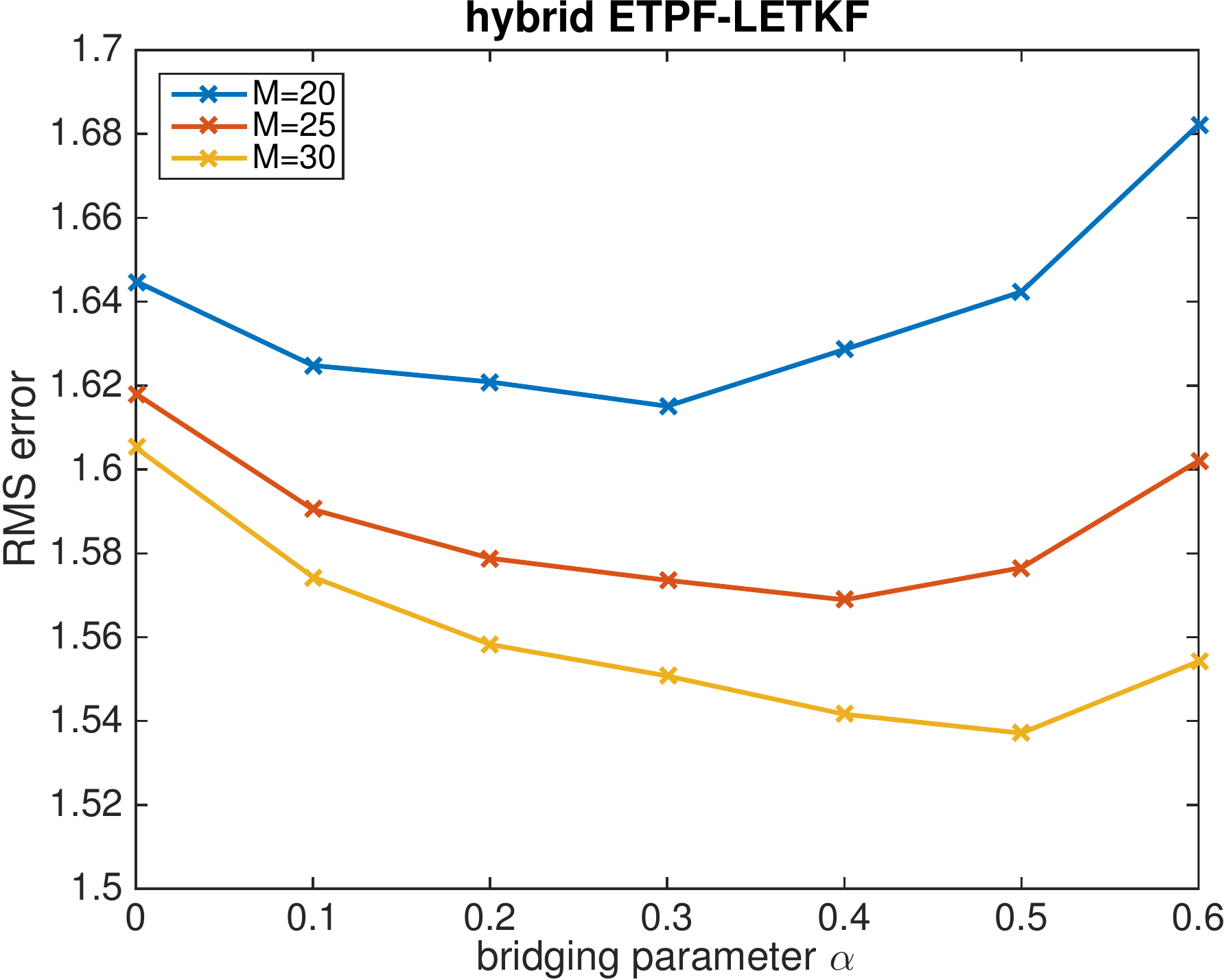} $\quad$
\includegraphics[width=0.45\textwidth]{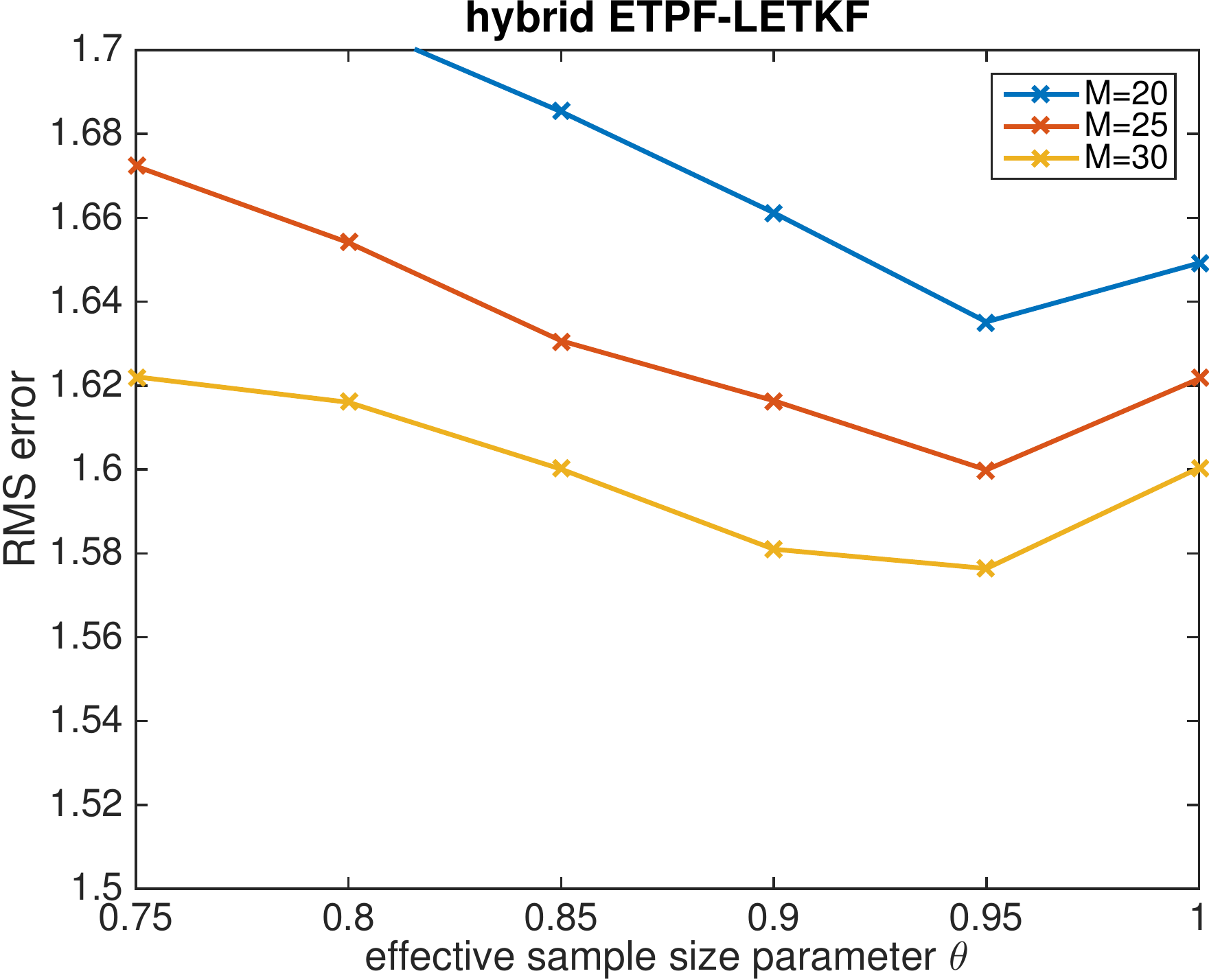} 
\end{center}
\caption{Hybrid filter for Lorenz-96 coupled to a wave equation with ensemble sizes $M=20,25,30$:
left panel: RMS errors for fixed bridging parameter; right panel: RMS errors 
for adaptively determined bridging parameters} \label{figure3}
\end{figure}

The  resulting time-averaged RMS errors are displayed in Figure \ref{figure3}. 
One can see that the minimal RMS error is 
obtained at $\alpha=0.3$ when $M=20$, at $\alpha=0.4$ when $M=25$ and at $\alpha=0.5$ when $M=30$.  The results for the adaptive approach show, that the RMS error is minimal 
for values of $\theta \in [0.90, 0.95]$ depending on the ensemble size, $M$, 
and that the adaptive approach is less efficient than an 
optimal choice of $\alpha$. This is likely due to the fact that an adaptive choice of the bridging parameter
in each grid point, $x_k$, together with localization leads to  less balanced fields. 

We finally mention that the damping in (\ref{wave}) is necessary to suppress the accumulation of imbalances generated in each assimilation step. See \cite{sr:br10b} for a detailed discussion 
of this phenomenon and a potential cure in the form of a mollified filter approach. 


\subsection{Convergence study for a single Bayesian inference step} \label{sec:convergence}

\begin{figure}
\begin{center}
\includegraphics[width=0.45\textwidth]{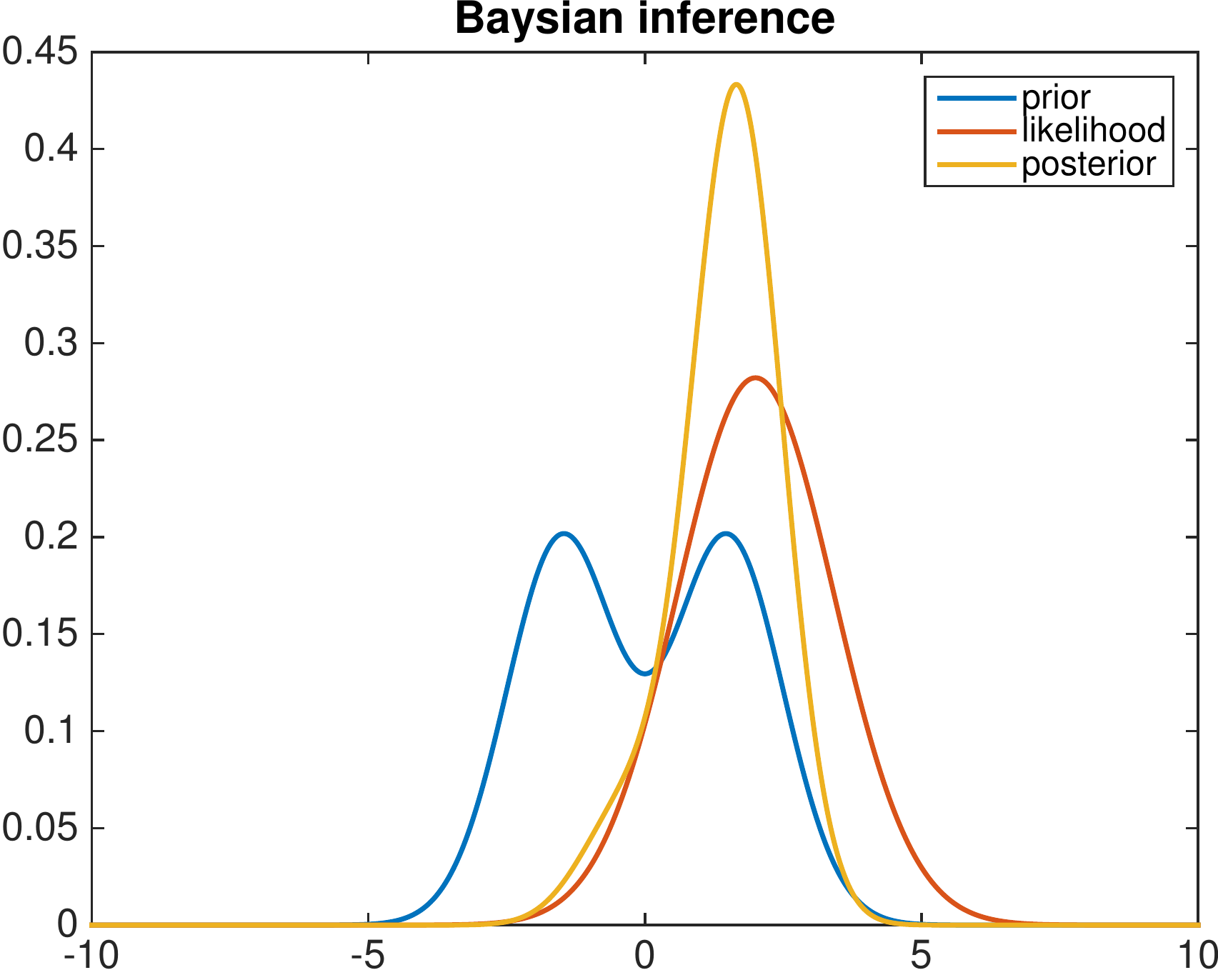} $\quad$
\includegraphics[width=0.45\textwidth]{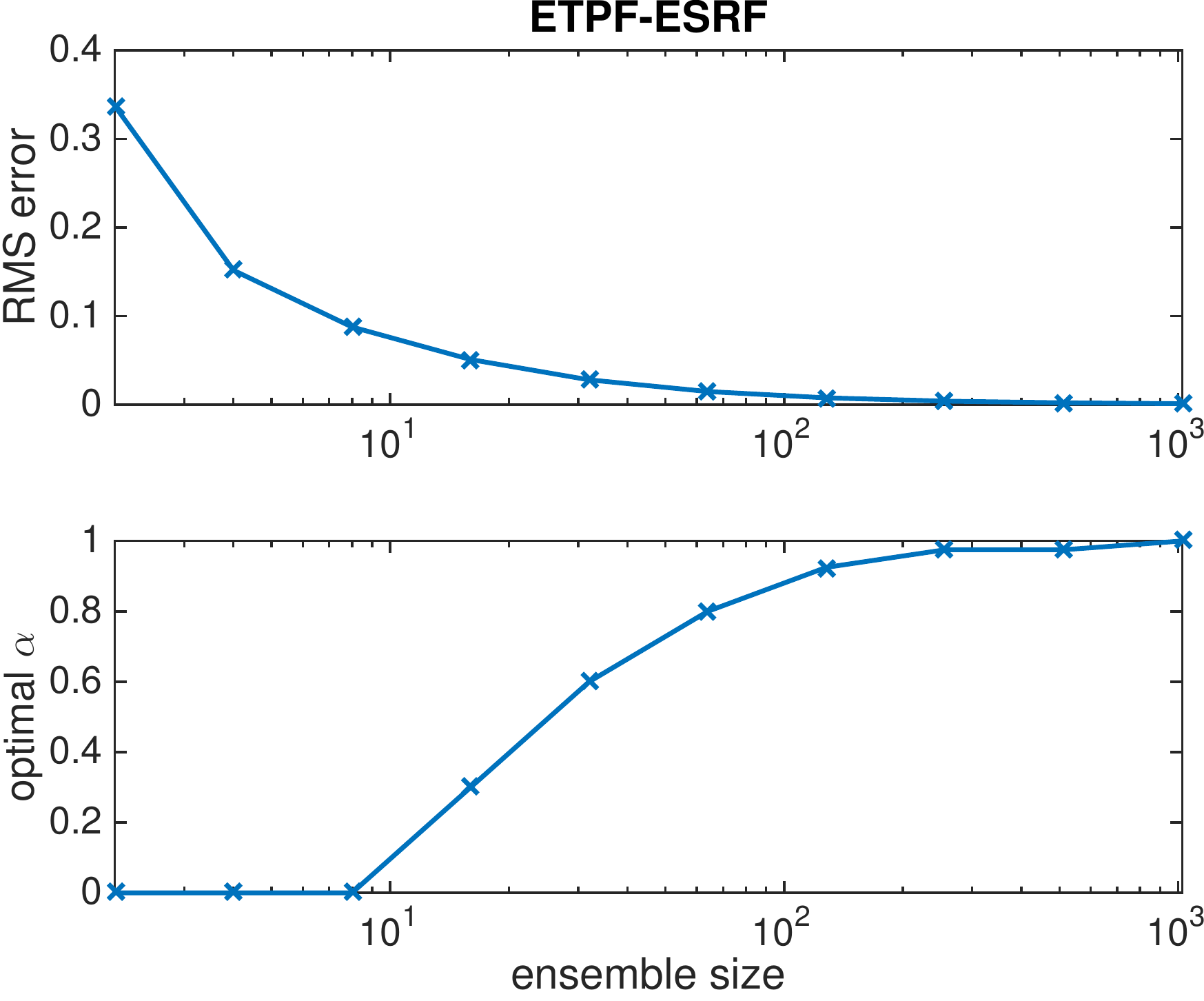} 
\end{center}
\caption{Single data assimilation step with ETPF-ESRF; left panel: prior and posterior distributions; right panel: optimal
bridging parameter and resulting RMS errors as a function of the ensemble size.} \label{figure4}
\end{figure}

After conducting a series of numerical experiments for which the ensemble sizes were relatively small, we finally
demonstrate that the proposed hybrid filter indeed bridges EnKFs and PFs as the ensemble size varies from very small
to very large. In order to keep the computations feasible, we only perform a singe assimilation step with given bimodal prior
(see Figure \ref{figure4}) and the ETPF-ESRF implemented for ensemble sizes $M \in [2,4,8,\ldots,1024]$. 
Each experiment is repeated 100,000 times
and the averaged RMS errors for the estimated posterior mean are computed in terms of $M$ and $\alpha$. The optimal
value of $\alpha$ and its RMS error are displayed in Figure \ref{figure4} as a function of the ensemble size, $M$. As expected,
the optimal value of $\alpha$ is close to zero for small $M$ and approaches one as $M$ increases. 

\section{Conclusions} Further numerical experiments show that the ETPF-ESRF and 
ETPF-LETKF, respectively, outperform the ESRF-ETPF and LETKF-ETPF, respectively, on all
numerical examples of Section \ref{sec:num}. See \cite{sr:nawinda} for more details. 
We believe that this is due to the fact that 
the true forecast/prior distributions in our examples are more non-Gaussian than the 
analysis/posterior 
distributions. This is the case, in particular, for the Lorenz 63 model. 
Under those circumstances it seems preferable to apply a PF first followed by an EnKF. 
We also found that
the adaptive hybrid method works well for the Lorenz-96 system, while a fixed choice 
of $\alpha$ is preferable for the Lorenz-63 model. It is feasible that other adaptation criteria for the choice of
$\alpha$ could lead to better results.  For example, one could select $\alpha$ such that the 
Kullback-Leibler divergence between the uniform measure, $1/M$, and the importance 
weights (\ref{weights}) remains below a given threshold. 

Our numerical results also indicate that our hybrid method is suitable for multi-scale systems. 
However, our hybrid method does not address the problem
of generating unbalanced fields through localization (see, for example, \cite{sr:kalnay09a}). 
This issue will be addressed in a separate
publication using the idea of mollification, as proposed in \cite{sr:br10b}.
We finally mention that our hybrid methods can be combined with modified proposal densities for 
PFs (see \cite{sr:leeuwen15} for a recent review of the subject).

\section*{Acknowledgments} SR and MR acknowledge support under the DFG Collaborative
Research Center SFB 1114: Scaling Cascades in Complex Systems. NC's work has been
supported by a fellowship of the Helmholtz Graduate School GeoSim. We thank Walter Acevedo
for discussions on the implementation of the hybrid method.


\bibliographystyle{plain}
\bibliography{survey}

\begin{thebibliography}{10}

\bibitem{sr:anderson10}
J.L. Anderson.
\newblock A non-{G}aussian ensemble filter update for data assimilation.
\newblock {\em Monthly Weather Review}, 138:4186--4198, 2010.

\bibitem{sr:crisan}
A.~Bain and D.~Crisan.
\newblock {\em Fundamentals of stochastic filtering}, volume~60 of {\em
  Stochastic modelling and applied probability}.
\newblock Springer-Verlag, New-York, 2008.

\bibitem{sr:kalnay09a}
J.~Ballabrera-Poy, E.~Kalnay, and S.-C. Yang.
\newblock Data assimilation in a system with two scales--combining two
  initialization techniques.
\newblock {\em Tellus}, 61A:539--549, 2009.

\bibitem{sr:bengtsson08}
T.~Bengtsson, P.~Bickel, and B.~Li.
\newblock Curse of dimensionality revisited: {C}ollapse of the particle filter
  in very large scale systems.
\newblock In {\em IMS Lecture Notes - Monograph Series in Probability and
  Statistics: Essays in Honor of David F.~Freedman}, volume~2, pages 316--334.
  Institute of Mathematical Sciences, 2008.

\bibitem{sr:br10b}
K.~Bergemann and S.~Reich.
\newblock A mollified ensemble {K}alman filter.
\newblock {\em Q. J. R. Meteorological Soc.}, 136:1636--1643, 2010.

\bibitem{sr:reich15}
Y.~Chen and S.~Reich.
\newblock Assimilating data into scientific models: {A}n optimal coupling
  perspective.
\newblock In {\em Frontiers in Applied Dynamical Systems: Reviews and
  Tutorials}, volume~2, pages 75--118. Springer-Verlag, New York, 2015.

\bibitem{sr:nawinda}
Nawinda Chutstagulprom.
\newblock Phd thesis.
\newblock Technical report, Universit\"at Potsdam, Institut f\"ur Mathematik,
  Am Neuen Palais 10, 14469 Potsdam, 2016.

\bibitem{sr:handbook}
D.~Crisan and B.~Rozovskii, editors.
\newblock {\em The Oxford Handbook of Nonlinear Filtering}.
\newblock Oxford University Press, Oxford, 2011.

\bibitem{sr:Doucet}
A.~Doucet, N.~de~Freitas, and N.~Gordon (eds.).
\newblock {\em Sequential Monte Carlo methods in practice}.
\newblock Springer-Verlag, Berlin Heidelberg New York, 2001.

\bibitem{sr:evensen}
G.~Evensen.
\newblock {\em Data assimilation. {T}he ensemble Kalman filter}.
\newblock Springer-Verlag, New York, 2006.

\bibitem{sr:frei11}
M.~Frei and H.R. K\"unsch.
\newblock Bridging the ensemble {K}alman and particle filters.
\newblock {\em Biometrica}, 100:781--800, 2013.

\bibitem{sr:gaspari99}
G.~Gaspari and S.E. Cohn.
\newblock Construction of correlation functions in two and three dimensions.
\newblock {\em Q. J. Royal Meteorological Soc.}, 125:723--757, 1999.

\bibitem{sr:hunt07}
B.R. Hunt, E.J. Kostelich, and I.~Szunyogh.
\newblock Efficient data assimilation for spatialtemporal chaos: {A} local
  ensemble transform {K}alman filter.
\newblock {\em Physica D}, 230:112--137, 2007.

\bibitem{sr:jazwinski}
A.H. Jazwinski.
\newblock {\em Stochastic processes and filtering theory}.
\newblock Academic Press, New York, 1970.

\bibitem{sr:kalnay}
E.~Kalnay.
\newblock {\em Atmospheric modeling, data assimilation and predictability}.
\newblock Cambridge University Press, 2002.

\bibitem{sr:stuart15}
K.~Law, A.~Stuart, and K.~Zygalakis.
\newblock {\em Data Assimilation}.
\newblock Springer-Verlag, New York, 2015.

\bibitem{sr:lei11}
J.~Lei and P.~Bickel.
\newblock A moment matching ensemble filter for nonlinear and non-{G}aussian
  data assimilation.
\newblock {\em Mon.~Weath.~Rev.}, 139:3964--3973, 2011.

\bibitem{sr:lorenz63}
E.N. Lorenz.
\newblock Deterministic non-periodic flows.
\newblock {\em J. Atmos. Sci.}, 20:130--141, 1963.

\bibitem{sr:lorenz96}
E.N. Lorenz.
\newblock Predictibility: {A} problem partly solved.
\newblock In {\em Proc. Seminar on Predictibility}, volume~1, pages 1--18,
  ECMWF, Reading, Berkshire, UK, 1996.

\bibitem{sr:lorenz98}
E.N. Lorenz and K.E. Emanuel.
\newblock Optimal sites for suplementary weather observations: {S}imulations
  with a small model.
\newblock {\em J. Atmos. Sci.}, 55:399--414, 1998.

\bibitem{sr:snyder13}
S.~Metref, E.~Cosme, C.~Snyder, and P.~Brasseur.
\newblock A non-{G}aussian analysis scheme using rank histograms for ensemble
  data assimilation.
\newblock {\em Nonlinear Processes in Geophysics}, 21:869--885, 2013.

\bibitem{sr:wright99}
J.~Nocedal and S.J. Wright.
\newblock {\em Numerical Optimization}.
\newblock Springer-Verlag, New York, 2nd edition, 2006.

\bibitem{sr:Pele-iccv2009}
O.~Pele and M.~Werman.
\newblock Fast and robust earth mover's distances.
\newblock In {\em Computer Vision, 2009 IEEE 12th international conference},
  pages 460--467, 2009.

\bibitem{sr:vanhandel15}
P.~Rebeschini and R.~{v}an Handel.
\newblock Can local particle filters beat the curse of dimensionality?
\newblock {\em Ann. Appl. Probab.}, 25:2809--2866, 2015.

\bibitem{sr:reich11}
S.~Reich.
\newblock A {G}aussian mixture ensemble transform filter.
\newblock {\em Q.~J.~R.~Meterolog.~Soc.}, 138:222--233, 2012.

\bibitem{sr:reich13}
S.~Reich.
\newblock A nonparametric ensemble transform method for {B}ayesian inference.
\newblock {\em SIAM J. Sci. Comput.}, 35:A2013--A2024, 2013.

\bibitem{sr:reichcotter15}
S.~Reich and C.J. Cotter.
\newblock {\em Probabilistic Forecasting and Bayesian Data Assimilation}.
\newblock Cambridge University Press, Cambridge, 2015.

\bibitem{sr:smith07}
K.W. Smith.
\newblock Cluster ensemble {K}alman filter.
\newblock {\em Tellus}, 59A:749--757, 2007.

\bibitem{sr:stordal11}
A.S. Stordal, H.A. Karlsen, G.~N\ae{}vdal, H.J. Skaug, and B.~Vall\'es.
\newblock Bridging the ensemble {K}alman filter and particle filters: the
  adaptive {G}aussian mixture filter.
\newblock {\em Comput. Geosci.}, 15:293--305, 2011.

\bibitem{sr:tippett03}
M.K. Tippett, J.L. Anderson, G.H. Bishop, T.M. Hamill, and J.S. Whitaker.
\newblock Ensemble square root filters.
\newblock {\em Mon. Wea. Rev.}, 131:1485--1490, 2003.

\bibitem{sr:leeuwen15}
P.J. {V}an Leeuwen.
\newblock Nonlinear data assimilation for high-dimensional systems.
\newblock In {\em Frontiers in Applied Dynamical Systems: Reviews and
  Tutorials}, volume~2, pages 1--73. Springer-Verlag, New York, 2015.

\end{thebibliography}

\end{document}